\documentclass[12pt]{amsart}
\usepackage{amssymb}
\usepackage{amscd}
\def\Aff{\mathop{\mathrm {Aff}}\nolimits}
\def\aff{\mathop{\mathrm {aff}}\nolimits}

\def\Ad{\mathop{\mathrm {Ad}}\nolimits}
\def\ad{\mathop{\mathrm {ad}}\nolimits}

\def\Lie{\mathop{\mathrm {Lie}}\nolimits}

\def\MD{\mathop{\mathrm {MD}}\nolimits}
\def\G{\mathop{\mathrm {G}}\nolimits}
\newtheorem{theorem}{Theorem}[subsection]
\newtheorem{proposition}[theorem]{Proposition}
\newtheorem{lemma}[theorem]{Lemma}

\newtheorem{corollary}[theorem]{Corollary}
\newtheorem{definition}[theorem]{DEFINITION}
\begin{document}
\title{Quantum co-adjoint orbits of $\MD_4$-groups.}
\author{ Nguyen Viet Hai}
\date {Version of March 10, 2000}
\address{Haiphong Teacher's Training College, Haiphong city, Vietnam}
\email{nguyen\_viet\_hai@yahoo.com}
\thanks{This work was supported in part by the Vietnam National  Foundation for Fundamental Science Research; \hspace{10.5cm}Typeset {\LaTeX}}

\keywords{Moyal $\star$-product,  $\MD_4$-group,  quantum half-plans,  quantum  rotation paraboloids, quantum hyperbolic cylinders,  quantum hyperbolic paraboloids.}
\maketitle
\begin{abstract}
Using $\star$-product on Co-adjoint orbits (K-orbits) of  the $\MD_4$- groups  we obtain quantum half-planes, quantum hyperbolic cylinders, quantum hyperbolic paraboloids...via Fedosov deformation quantization.  From this we have corresponding unitary representations of  the  $\MD_4$- groups. Particularly, for groups $\G_{4,2,3(\varphi)};  \G_{4,2,4};  \G_{4,3,4(\varphi)}$ and $ \G_{4,4,1} $, which are neither nilpotent nor exponential ,  we  obtain the explicit formulas.   
\end{abstract}
\section{Introduction}
 First of all,  we recall the notion of K-action (see \cite{kirillov1}).  Let us  denote  by $\G$ a connected and simply connected Lie group,   its Lie algebra  $\mathfrak g=T_e\G$ as the tangent space at the neutral element e.  It is easy to see that to each element $g \in \G  $ one can associate a map $$A(g):  \G \longrightarrow \G  $$  by conjugacy,  in fixing the neutral element $e \in \G$. Therefore, there is a corresponding tangent map $$A(g)_*:  \mathfrak g=T_e\G \longrightarrow \mathfrak g=T_e\G$$
$$X\in \mathfrak g \mapsto \frac{d}{dt}g\exp(tX)g^{-1}\vert_{t=0} \in\mathfrak g. $$
 This  defines an action,  denoted as usually by $\Ad$, of group  $\G$  in  its Lie algebra $\mathfrak g$. We define the  co-adjoint action of group $\G$ in the dual vector space $\mathfrak g^*$ by the formula: 
$$\langle K(g)F,X \rangle:=\langle F,Ad(g^{-1})X \rangle $$
for all $F\in \mathfrak g^*, X\in \mathfrak g , g\in \G$.  It is easy to check that this defines an action of  $\G$ on $\mathfrak g^*.$
As an easy consequence, the dual space $\mathfrak g^*$ is decomposed into a disconnected sum of the K-orbits. In 1980  Do Ngoc Diep  introduced the notion of  the $\MD$-groups, $\MD$-algebras and then Le Anh Vu gave a complete classification of the  $\MD_4$-groups (see \cite{dndiep}). In  \cite {diephai1}, \cite{diephai2}  and   \cite{hai3} we  obtained  the exact formulae of deformation quantization and therefore performed the corresponding quantun coadjoint orbits for the $\overline{\MD}$-groups  (i.e. the groups,  every K-orbit of which is of dimension, equal 0 or  dim $\G$) and for the real  diamond Lie group.   In this article,  applying the procedure of deformation quantization we shall obtain  quantum co-adjoint orbits of all  the $\MD_4$-groups.\par
The article is organized as follows.  In section 2 we   recall the basic definitions, preliminary results. Each adapted chart that carries the  Moyal $\star$-product from $\mathbb R^2$  onto co-adjoint orbits $\Omega_F$ of the exponential $\MD_4$-groups,  in particular,  Hamiltonian functions in canonical coordinates,  are introduced in section 3.  Section 4 is  devoted to the following  groups $$\G_{4,2,3(\frac{\pi}{2})}; \G_{4,2,4}; \G_{4,3,4(\frac{\pi}{2})}; \G_{4,4,1}.$$    By direct  computations  and  by exponentiating we obtain the corresponding unitary representations of the  $\MD_4$-groups.
\section{Basic definitions and Preliminary results}
\begin{definition}(\cite{dndiep}) We say that a solvable Lie group $\G$  belongs to the  class $ \MD$ if and only if  every its K-orbit has dimension 0 or maximal. A Lie algebra is of class   $ \MD$ if and only if its corresponding Lie group is of the same class.\end{definition}
We also recall  the following results for   $\MD_4$-algebras   (i.e. $  dim\  \mathfrak   g= 4 $ ):
\begin{theorem}
Assume $\mathfrak g$ is a  $ \MD_4$-algebra with generators $  X,Y,Z,T$,  then\par
\begin{itemize}
\item[I.] If $\mathfrak g$ is decomposable  then it is of the form 
$$\mathfrak g=\mathbb R^n\oplus\tilde{\mathfrak g}$$ for  n=1,2,3,4 and some indecomposable ideal $\tilde{\mathfrak g}.$\par
\item[II.] If $\mathfrak g$ is indecomposable then $ \mathfrak g$ is of class $ \MD_4$ if and only if it is generated by the generators $ X,Y,Z,T$  with only non-trivial commutation relation which is one of following relations defined in each case: 
\end{itemize}
\begin {itemize}
\item[1.]  $\mathfrak g^1=[\mathfrak g,\mathfrak g]=\mathbb RZ\cong\mathbb R,$  and

   \begin{itemize}
        \item[1.1] $ [T,X]=Z$ \hspace{10cm} $(\mathfrak g_{4,1,1})$
        \item[1.2] $ [T,Z]=Z$ \hspace{10cm} $(\mathfrak g_{4,1,2})$
   \end{itemize}
\item[2.]  $\mathfrak g^1=[\mathfrak g,\mathfrak g]=\mathbb RY+\mathbb RZ\cong\mathbb R^2$, and
   \begin{itemize}
       \item[2.1] $[T,Y]=\lambda Y,[T,Z]=Z;\lambda \in \mathbb R^*=\mathbb R\backslash (0)\hspace{4.8cm} (\mathfrak g_{4,2,1(\lambda)})$
       \item[2.2] $[T,Y]=Y;[T,Z]=Y+Z   \hspace{7.5cm} (\mathfrak g_{4,2,2})$
       \item[2.3] $\ad_T=\left(\begin{array}{ccc} cos\varphi & sin\varphi & 0 \\ -sin\varphi & cos\varphi & 0\\0 & 0 &0 \end{array} \right)$,\hspace{6.3cm} $(\mathfrak g_{4,2,3(\varphi)})$
      \item[2.4] $\ad_T=\left(\begin{array}{ccc}1 & 0 & 0 \\ 0 & 1 &  0\\ 0 & 0 & 0 \end{array} \right),$
               $\ad_X=\left(\begin{array}{ccc} 0 & 1 & 0 \\ -1 & 0  &0 \\ 0 &0 & 0\end{array} \right)$                                             \hspace{1.5cm} $(\mathfrak g_{4,2,4})=\Lie(\Aff(\mathbb C))$
    \end{itemize}
\item[3.] $\mathfrak g^1=[\mathfrak g,\mathfrak g]=\mathbb RX+\mathbb RY+\mathbb RZ$, and 
    \begin{itemize}
        \item[3.1] $\ad_T=\left(\begin{array}{ccc} \lambda_1  & 0 & 0 \\ 0 & \lambda_2  & 0 \\ 0 &0 & 1\end{array} \right)$, $\lambda_1,\lambda_2 \in \mathbb R^*$   \hspace{4.7cm} $(\mathfrak g_{4,3,1(\lambda_1,\lambda_2)})$
        \item[3.2] $\ad_T=\left(\begin{array}{ccc} \lambda  & 1 & 0 \\ 0 & \lambda  & 0 \\ 0 &0 & 1\end{array} \right)$,$\lambda \in \mathbb R^*$  \hspace{6.5cm} $(\mathfrak g_{4,3,2(\lambda)})$
        \item[3.3] $\ad_T=\left(\begin{array}{ccc} 1  & 1 & 0 \\ 0 & 1  & 1 \\ 0 &0 & 1\end{array} \right)$\hspace{8.5cm} $(\mathfrak g_{4,3,3})$
        \item[3.4]  $\ad_T=\left(\begin{array}{ccc} cos\varphi & sin\varphi & 0 \\ -sin\varphi & cos\varphi & 0\\0 & 0 & \lambda \end{array} \right)$,$\lambda\in \mathbb R^*,\varphi\in (0,\pi)$                           \hspace{3cm} $(\mathfrak g_{4,3,4(\lambda)})$

     \end{itemize}
\item[4.]  $\mathfrak g^1=[\mathfrak g,\mathfrak g]=\mathbb RX+\mathbb RY+\mathbb RZ\cong \mathfrak h_3$  -the 3-dimensional Heisenberg Lie algebra,  and 
    \begin{itemize}
       \item[4.1]  $\ad_T=\left(\begin{array}{ccc}0  & 1 & 0 \\ -1 & 0  & 0 \\ 0 &0 & 0\end{array} \right)$, $[X,Y]=Z$                  \hspace{3cm} $(\mathfrak g_{4,4,1}=\Lie(\mathbb R \ltimes_j\mathbb H_3))$
       \item[4.2] $\ad_T=\left(\begin{array}{ccc}-1  & 0 & 0 \\ 0 & 1  & 0 \\ 0 &0 & 0\end{array} \right)$, $[X,Y]=Z$                  \hspace{3cm} $(\mathfrak g_{4,4,2}=\Lie(\mathbb R \ltimes\mathbb H_3))$\par
(In this case the group is called the real diamond Lie  group, see\cite{hai3})
     \end{itemize}
\end{itemize} 
\end{theorem}
Until now we fix a  $ \MD_4$- algebra $\mathfrak g$ with the standard basis X,Y,Z,T. It is  isomorphic to $\mathbb R^4$ as vector spaces. The coordinates in this standard basis is denoted by (a,b,c,d). We identify its dual vector spase $\mathfrak g^*$ with $\mathbb R^4$ with the help of the dual  basis  $X^*,Y^*,Z^*,T^*$  and with the local coordinates $(\alpha,\beta,\gamma,\delta)$. Thus, for  all  $ U\in \mathfrak g, U=aX+bY+cZ+dT$ and for all $ F\in\mathfrak g^*,F=\alpha X^*+\beta Y^*+\gamma Z^*+\delta T^*.$ Finally, $\Omega_F$  is the co-adjoint  orbit passing through $F\in\mathfrak g^*.$
\begin{theorem}(The Picture of Co-adjoint Orbit)\cite{dndiep}
\begin {itemize}
\item[1.1] Case\quad $ \G=\G_{4,1,1}$
   \begin{itemize}
        \item[i.] Each point F with the coordinate  $\gamma=0$ is a  0-dimensional co-adjoint orbit:          $$\Omega_F=\Omega_{(\alpha,\beta,0,\delta)}$$
        \item[ii.] The subset $\gamma \ne 0 $ is decomposed into  a family of  2-dimensional co-adjoint orbits:               \begin{equation}\Omega_F=\Omega_{(\beta,\gamma \ne 0)}=\{(\alpha+\gamma d,\beta,\gamma,-\gamma a+\delta)\}=\{(x,\beta,\gamma,t)\vert x,t\in\mathbb R\},\end{equation}
which are{\quad\bf planes.}
   \end{itemize}
\item[1.2.] Case\quad $ \G=\G_{4,1,2}$
   \begin{itemize}
       \item[i.]  Each point F with the coordinate  $\gamma=0$ is a  0-dimensional co-adjoint orbit:          $$\Omega_F=\Omega_{(\alpha,\beta,0,\delta)}$$
       \item[ii.] The subset $\gamma \ne 0$ is decomposed into a family  of 2-dimensional co-adjoint orbits:
\begin{equation} \Omega_F=\Omega_{(\alpha,\beta)}=\{\alpha,\beta ,\gamma e^d,-\gamma c\sum_1^{\infty}\frac{d^{n-1}}{n!}+\delta\}=\{(\alpha,\beta,z,t)\vert z,t\in\mathbb R,\gamma z>0 \}\end{equation} which are { \bf  half-planes}, parameterized by the coordinates $\alpha,\beta \in \mathbb R$ 
    \end{itemize}
\item[2.1] Case\quad $\G=\G_{4,2,1(\lambda)},\lambda \in \mathbb R^*$
    \begin{itemize}
        \item[i.]  Each point on the plane \quad$\beta=\gamma=0$ is a  0-dimensional co-adjoint orbit:          $$\Omega_F=\Omega_{(\alpha,0,0,\delta)}$$
        \item[ii.]  The open set\quad $\beta^2+\gamma^2 \ne 0$ is decomposed into the union of 2-dimensional{ \bf cylinders} 
\begin{equation}\Omega_F=\{(\alpha,\beta e^{s\lambda},\gamma e^s,t) \vert s,t \in\mathbb R \}\end{equation}
     \end{itemize}
\item[2.2] Case\quad $\G=\G_{4,2,2}$  
    \begin{itemize}
       \item[i.]  Each point on the plane\quad $\beta=\gamma=0$ is a  0-dimensional co-adjoint orbit:          $$\Omega_F=\Omega_{(\alpha,0,0,\delta)}$$
       \item[ii.] The open set \quad$\beta^2+\gamma^2 \ne 0$ is decomposed into the union of 2-dimensional{ \bf cylinders}: 
\begin{equation}\Omega_F=\{(\alpha,\beta e^s,\beta se^s+\gamma e^s,t) \vert s,t \in\mathbb R \}\end{equation}
     \end{itemize}
\item[2.3]  Case\quad $\G=\G_{4,2,3(\varphi)}$  with  $\varphi \in(0,\pi). $  We identify $\mathfrak g_{4,2,3(\varphi)}^*$ with $\mathbb R\times \mathbb C\times\mathbb R$ and $F=(\alpha,\beta,\gamma,\delta )$ with  $(\alpha,\beta+i\gamma,\delta).$ Then,
    \begin{itemize}
       \item[i.] Each point $(\alpha,0,\delta)$ is a  0-dimensional co-adjoint orbit:  $$\Omega_F=\Omega_{(\alpha,0+i0,\delta)}$$
       \item[ii.]  The open set \quad $\beta +i\gamma \ne 0$  is decomposed into the union of 2-dimensional
  co-adjoint orbits:
\begin{equation} \Omega_F=\{(\alpha,(\beta+i\gamma)e^{se^{i\varphi}},t)\vert s,t \in\mathbb R\},\end{equation} which are also{ \bf cylinders.}
   \end{itemize}
 \item[2.4] Case\quad $\G=\G_{4,2,4}=\widetilde\Aff(\mathbb C)$
     \begin{itemize}
     \item[i.]  Each point $(\alpha,0,0,\delta)$ is a  0-dimensional co-adjoint orbit: $$\Omega_F =  \Omega_{(\alpha,0,0,\delta)}$$
      \item[ii.] The open set $\beta^2+\gamma^2 \ne 0$ is the single 4-dimensional co-adjoint orbit: \begin{equation}\Omega_F=\Omega_{(\beta^2+\gamma^2) \ne 0}=\{(x,y,z,t)\vert y^2+z^2 \ne 0\}=\mathbb R\times(\mathbb R^2)^*\times\mathbb R\end{equation}
      \end{itemize}
\item[3.1]  Case $\G$ is one of the groups \quad $\G_{4,3,1(\lambda_1,\lambda_2)},$   $\G_{4,3,2(\lambda)}$  or\quad  $\G_{4,3,3}$
   \begin{itemize}
      \item[i.] Each point  $F=\delta T^*$ on the line   $\alpha=\beta=\gamma=0$ is a 0-dimensional co-adjoint orbit.
      \item[ii.] The open set\quad  $\alpha^2+\beta^2+\gamma^2 \ne 0$  is decomposed into a family of  co-adjoint orbits which are { \bf cylinders},  coressponding to the groups\quad  $\G_{4,3,1(\lambda_1,\lambda_2)}$,$\G_{4,3,2(\lambda)}$  or  $\G_{4,3,3}$:
 \begin{equation}\Omega_F =\{(\alpha e^{s\lambda_1},\beta e^{s\lambda_2},\gamma e^s,t)\vert s,t \in\mathbb R\}.\end{equation}
 \begin{equation}\Omega_F=\{(\alpha e^{s\lambda},\alpha se^{s\lambda}+\beta e^{s\lambda},\gamma e^s,t)\vert s,t\in\mathbb R\}.\end{equation}
 \begin{equation}\Omega_F=\{(\alpha e^{s\lambda},\alpha se^s+\beta e^s,\frac 12 \alpha s^2e^s+\beta se^s+\gamma e^s,t)\vert s,t\in\mathbb R\}.\end{equation}
    \end{itemize}
\item[3.4]  Case\quad $\G=\G_{4,3,4(\lambda,\varphi)}$  for $\lambda \in \mathbb R^*$,  $\varphi\in(0,\pi).$  We identify $\mathfrak g^*_{4,3,4(\lambda,\varphi)}$ with $\mathbb C \times \mathbb R^2$ and  $F=(\alpha,\beta,\gamma,\delta)$ with $(\alpha+i\beta,\gamma,\delta).$ Then,
   \begin{itemize}
      \item[i.]  Each point of the line defined by the condition $\alpha=\beta=\gamma=0$ is a  0-dimensional co-adjoint orbit  $$\Omega_F= \Omega_{(0,0,\delta)}=\{(0+i.0,0,\delta)\}.$$

       \item[ii.]  The open set $\vert \alpha+i\beta\vert^2+\gamma^2 \ne 0$ is decomposed into an union of co-adjoint orbits, which are  { \bf cylinders}
   \begin{equation}\Omega_F=\{((\alpha+i\beta)e^{se^{i\varphi}},\gamma e^{s\lambda},t)\vert s,t\in\mathbb R \}.\end{equation}
    \end{itemize} 
\item[4.1] Case\quad $\G=\G_{4,4,1}=\mathbb R\ltimes_j\mathfrak h_3$
   \begin{itemize}
      \item[i.]  Each point of the line defined by the conditions $\alpha=\beta=\gamma=0$  is a  0-dimensional      co-adjoint orbit
 $$\Omega_F= \Omega_{(0,0,0,\delta)}=\{(0,0,0,\delta)\}.$$
      \item[ii.] The set $\alpha^2+\beta^2 \ne 0,\gamma=0 $  is the union of 2- dimensional co-adjoint orbits, which are { \bf rotation cylinders}
\begin{equation}\Omega_F=\{(\alpha cos\theta-\beta sin\theta ,\alpha sin\theta+\beta cos\theta,0,t)\vert \theta,t\in \mathbb R \}\end{equation} i.e  $$\Omega_F=\{(x,y,0,t)\vert x^2+y^2=\alpha^2+\beta^2; \quad x,y,t\in\mathbb R \}.$$
      \item[iii.]  The open set $\gamma \ne 0$ is decomposed into a  union of 2-dimensional co-adjoint orbits 
\begin{equation}\Omega_F= \{(x,y,\gamma,t)\vert x^2+y^2-2\gamma t=\alpha^2+\beta^2-2\gamma\delta ;\quad  x,y,t\in\mathbb R \},\end{equation} which are { \bf rotation paraboloids}
   \end{itemize}
\item[4.2]  Case\quad $\G=\G_{4,4,2}=\mathbb R \ltimes \mathbb H_3$,  the real diamond group (see  \cite {hai3})
 \begin{itemize}
     \item[i.] Each point of the line $\alpha=\beta=\gamma =0$ is a  0-dimensional co-adjoint orbit
$$\Omega_F=\Omega_{(0,0,0,\delta)}$$
     \item[ii.] The set $\alpha\ne 0,\beta=\gamma=0$ is union of  2-dimensional co-adjoint orbits ,which are just { \bf half-planes }
\begin{equation}\Omega_F=\{(x,0,0,t)\quad\vert\quad x,t\in \mathbb R , \alpha x>0\}\end{equation}
     \item[iii.]  The set\quad $\alpha=\gamma=0,\beta\ne 0$  is union of 2-dimensional co-adjoint orbits, which are just {  \bf half-planes}
\begin{equation}\Omega_F=\{(0,y,0,t)\quad\vert\quad y,t\in\mathbb R, \beta y >0 \}.\end{equation}
    \item[iv.] The set $\alpha\beta\ne 0,\gamma=0$ is decomposed into a family of  2-dimensional co-adjoint orbits, which are just{ \bf hyperbolic-cylinders}
\begin{equation}\Omega_F=\{(x,y,0,t)\quad\vert x,y,t\in\mathbb R\quad\&\quad \alpha x>0,\beta y>0,xy=\alpha\beta\}.\end{equation}
    \item[v.]  The open set  $ \gamma\ne 0$ is decomposed into a family of  2-dimensional co-adjoint orbits , which are just {  \bf  hyperbolic- paraboloids}
\begin{equation}\Omega_F=\{(x,y,\gamma,t)\quad\vert x,y,t \in\mathbb R \quad\&\quad xy- \alpha\beta=\gamma(t-\delta)\}.\end{equation}
    \end{itemize}     
   \end{itemize}
\end{theorem}
Thus,  we have 15 family of 2-dimensional co-adjoint orbits  and a  4-dimensional co-adjoint orbit $\Omega_{F} \cong {\mathbb C} \times {{\mathbb C}}^*$. They are strictly homogeneous symplectic manifolds with a flat action.\par
Deformation of Poisson brackets and associative algebras of $C^{\infty}$-functions on symplectic manifold (classical phase spaces) permit an autonomous quantization theory without need for Hilbert space operators. This theory is based on the of $\star$-products introduced by Flato and Lichnerowicz (see for instance \cite{arnalcortet1}).  Let us denote by $\Lambda$ the 2-tensor associated with the Kirillov standard form $\omega = dp \wedge dq$ in canonical Darboux coordinates. We  consider the well-known Moyal $\star$-product of two smooth functions $u,v \in C^{\infty}(\mathbb R^2)$  ( see e.g \cite{diephai1},  \cite{hai3}), defined by$$u \star v = u.v + \sum_{r \geq 1} \frac{1}{r!}(\frac{1}{2i})^r P^r(u,v),$$  where
$$P^1(u,v)=\{u,v\}$$
$$P^r(u,v) := \Lambda^{i_1j_1}\Lambda^{i_2j_2}\dots \Lambda^{i_rj_r}\partial^r_{i_1i_2\dots i_r} u \partial^r_{j_1j_2\dots j_r}v,$$ with $$\partial^r_{i_1i_2\dots i_r} := \frac{\partial^r}{\partial x^{i_1}\dots \partial x^{i_r}};\quad x:= (p,q) = (p_1,\dots,p_n,q^1,\dots,q^n).$$  It is well-known that this series converges in the Schwartz distribution spaces $\mathcal S (\mathbb R^n)$.\par
Remark that the $\MD_4$-groups are not all  nilpotent or exponential groups. In the most general context,  some quantum co-adjoint orbits appeared in \cite{arnalcortet1}\cite{arnalcortet2}. However, it is difficult to calculate precisely the  $\star$-product in concrete cases. In this article we will give  explicit formulas for  co-adjoint orbits of  all  $\MD_4$-groups,  even for groups which are neither nilpotent nor exponential.\par
\section{Quantum co-adjoint orbits of the exponential $\MD_4$-groups.}
In this  all section, we denote by $\G$  one of the following groups: (with $\varphi\ne\frac{\pi}{2}$) \par
$\G_{4,1,1}; \G_{4,1,2}; \G_{4,2,1(\lambda)}; \G_{4,2,2}; \G_{4,2,3(\varphi)}; \G_{4,3,1(\lambda_1,\lambda_2)}; \G_{4,3,2(\lambda)}; \G_{4,3,3};\G_{4,3,4(\lambda,\varphi)}; \G_{4,4,2}.$

\subsection{Hamiltonian functions in canonical coordinates of $\Omega_F $}\par
Each element $A\in \mathfrak g$ can be considered as a linear function $\widetilde A$ on co-adjoint orbits $(\subset \mathfrak g^*)$:\quad  $\widetilde A(F'):=\langle F',A \rangle,\quad F'\in \Omega_F$.  It is well-known that this function is the Hamiltonian function associated with the Hamiltonian vector field $\xi_A$, is defined on $\Omega_F$ by 
$$(\xi_Af)(x):=\frac{d}{dt}f(x\exp(tA))\vert_{t=0},\forall f \in C^{\infty}(\Omega_F).$$
The Kirillov form $\omega_F$ is defined by the formula 
\begin{equation}
\omega_F(\xi_A,\xi_B)=\langle F,[A,B]\rangle,\forall A,B  \in\mathfrak g.
\end{equation}
Denote by $\psi$ the indicated symplectomorphism from $\mathbb R^2$ onto $\Omega_F.$
$$(p,q)\in \mathbb R^2 \mapsto\psi(p,q)\in\Omega_F.$$

\begin{proposition} Each nontrivial orbit $\Omega_F\subset\mathfrak g^*$ of the co-adjoint representation of $\G$, admits a global  diffeomorphism $ \psi$
$$ \psi:\qquad (p,q)\in \mathbb R^2 \longmapsto \psi(p,q)\in\Omega_F$$ such that: 
\begin{itemize}
\item[\bf i.]    Hamiltonian function $\widetilde A=\langle F', A\rangle,\quad (F'\in \Omega_F; A\in \mathfrak g)$  is of  the form:
$$\widetilde A\circ \psi(p,q)=  \Phi(q).p+\Psi(q), $$ where   $ \Phi(q),\Psi(q)$ are   $C^\infty$ -functions on $\mathbb R$.\par
\item[\bf ii.] The  Kirillov  form (17)   is
\begin{equation} \omega=dp\wedge dq \end{equation}
\end{itemize}
\end{proposition}
{\it Proof.} 
\begin{itemize}
 \item[\bf i.] The diffeomorphism $\psi$  will be chosen case by case.
\begin{enumerate}
\item   Case   $ \G_{4,1,1} $  and  $\Omega_F$   defined by (1).\par
We chose the diffeomorphism:
$$\psi:\qquad \mathbb R^2\longrightarrow\Omega_F;\quad(p,q)\longmapsto  ( q,\beta,\gamma,p).$$ Then,  
 \begin{equation}\widetilde A\circ \psi(p,q)=dp+(aq+b\beta+c\gamma).\end{equation}
\item  Case   $\G_{4,1,2}$ and $\Omega_F $    defined by   (2).Then we take \par
$$\psi:\qquad \mathbb R^2\longrightarrow\Omega_F;\quad(p,q)\longmapsto(\alpha,\beta,\gamma e^q,p),$$
 \begin{equation}\widetilde A\circ \psi(p,q)=dp+(c\gamma e^q +a\alpha+b\beta).\end{equation}
\item  Case    $\G_{4,2,1(\lambda)}$  and $\Omega_F$   defined by  (3).
$$\psi:\qquad \mathbb R^2\longrightarrow\Omega_F;\quad(p,q)\longmapsto(\alpha,\beta e^{q\lambda},\gamma e^q,p),$$
 \begin{equation}\widetilde A\circ \psi(p,q)=dp+(c\gamma e^q  +a\alpha+b\beta e^{q\lambda}).\end{equation}
\item  Case   $\G_{4,2,2}$ and $\Omega_F $   defined by (4).
$$\psi:\qquad \mathbb R^2\longrightarrow\Omega_F;\quad  (p,q)\longmapsto(\alpha,\beta e^q,\beta qe^q+\gamma e^q,p),$$
 \begin{equation}\widetilde A\circ \psi(p,q)=dp+c\beta qe^q+(b\beta+c\gamma)e^q+a\alpha.\end{equation}
\item  Case   $\G_{4,2,3(\varphi),\varphi\ne \frac{\pi}{2}}$  and  $\Omega_F$    defined by  (5).
$$\psi:\qquad \mathbb R^2\longrightarrow\Omega_F;\quad(p,q)\longmapsto(\alpha, (\beta+i\gamma)e^{qe^{i\varphi}},p),$$
 \begin{equation}\widetilde A\circ \psi(p,q)=dp+(b+ic)(\beta+i\gamma)e^{qe^{i\varphi}}+a\alpha. \end{equation}
\item  Case   $\G_{4,3,1(\lambda_1,\lambda_2)}$  and $\Omega_F$    defined by  (7).
$$\psi:\qquad \mathbb R^2\longrightarrow\Omega_F;\quad(p,q)\longmapsto(\alpha e^{q\lambda_1},\beta e^{q\lambda_2},\gamma e^q,p),$$
 \begin{equation}\widetilde A\circ \psi(p,q)=dp+a\alpha e^{q\lambda_1}+b\beta e^{q\lambda_2}+c\gamma e^q. \end{equation}
\item  Case   $\G_{4,3,2(\lambda)}$  and $\Omega_F$   defined by  (8).
$$\psi:\qquad \mathbb R^2\longrightarrow\Omega_F;\quad(p,q)\longmapsto(\alpha e^{q\lambda},\alpha qe^{q\lambda}+\beta e^{q\lambda},\gamma e^q,p),$$
 \begin{equation}\widetilde A\circ \psi(p,q)=dp+(a\alpha+bq\alpha+b\beta)e^{q\lambda}+c\gamma e^q.  \end{equation}
\item  Case   $\G_{4,3,3}$  and  $\Omega_F$   defined by  (9).
$$\psi:\qquad \mathbb R^2\longrightarrow\Omega_F;\quad(p,q)\longmapsto(\alpha e^q,\alpha qe^q+\beta e^q,\frac 12\alpha q^2e^q+\beta qe^q+\gamma e^q,p),$$
 \begin{equation}\widetilde A\circ \psi(p,q)=dp+(a\alpha+b\alpha+b\beta+\frac 12c\alpha q^2+c\beta q+c\gamma)e^q. \end{equation}
\item  Case   $\G_{4,3,4(\lambda,\varphi)},\varphi\ne\frac{\pi}{2}$  and  $\Omega_F$   defined by  (10).
$$\psi:\qquad \mathbb R^2\longrightarrow\Omega_F;\quad(p,q)\longmapsto     ((\alpha+i\beta)e^{qe^{i\varphi}},\gamma e^{q\lambda},p),$$
 \begin{equation}\widetilde A\circ \psi(p,q)=dp+(a+ib)(\alpha+i\beta)e^{qe^{i\varphi}}+c\gamma e^{q\lambda}. \end{equation}
\item   Case $\G_{4,4,2}=\mathbb R \ltimes \mathbb H_3$ and\par 
        10.1 $\Omega_F$  defined by  (13).$$(p,q) \in \mathbb R^2 \mapsto \psi(p,q)=(\alpha e^{-q},0,0,p) \in \Omega_F,$$
\begin{equation}\tilde{A}\circ\psi(p,q) = dp +a\alpha e^{-q}. \end{equation}
      10.2 $\Omega_F$   defined by  (14).$$(p,q) \in \mathbb R^2 \mapsto \psi(p,q)=(0,\beta e^q,0,p) \in \Omega_F ,$$
 \begin{equation}\tilde{A}\circ\psi(p,q) =  dp +b\beta e^q.\end{equation}
      10.3  $\Omega_F$    defined by   (15).$$(p,q) \in \mathbb R^2 \mapsto \psi(p,q)=(\alpha e^{-q},\beta e^{q},0,p) \in \Omega_F,$$
\begin{equation}\tilde{A}\circ\psi(p,q) = dp +a\alpha e^{-q}+b\beta e^q. \end{equation}
       10.4  $\Omega_F$    defined by   (16) .$$   (p,q) \in \mathbb R^2 \mapsto \psi(p,q)=(e^{-q},(\alpha\beta+\gamma p-\gamma\delta)e^q,\gamma,p) \in \Omega_F,$$ \begin{equation}\tilde{A}\circ\psi(p,q) =ae^{-q}+b(\alpha\beta +\gamma p-\gamma\delta)e^q+c\gamma + dp=\end{equation}
$$ =(d+b\gamma e^q)p+ ae^{-q}+b(\alpha\beta-\gamma\delta)e^q+c\gamma.$$
(see \cite{hai3}  for more detail ). 
  \end{enumerate}
\item[\bf ii.]We prove the Kirillov  form on $\Omega_F$ is $\omega=dp\wedge dq$ , namely for the  case     $\G_{4,2,3(\varphi)},\varphi\ne\frac{\pi}{2}$.\par
From the Hamiltonian function $\tilde A\circ\psi(p,q)$    we have $$\xi_A(f)=\{\tilde A,f\}=d\frac{\partial f}{\partial q}-(b+ic)(\beta+i\gamma)e^{i\varphi}e^{qe^{i\varphi}}\frac{\partial f}{\partial p},$$ with  $A=aX+bY+cZ+dT \in \mathfrak g_{4,2,3(\varphi)}$
$$\xi_B(f)=\{\tilde B,f\}=d'\frac{\partial f}{\partial q}-(b'+ic')(\beta+i\gamma)e^{i\varphi}e^{qe^{i\varphi}}\frac{\partial f}{\partial p},$$  with $B=a'X+b'Y+c'Z+d'T \in \mathfrak g_{4,2,3(\varphi)}$.\quad Thus,
$$\xi_A\otimes\xi_B=dd'\frac{\partial}{\partial q}\otimes\frac{\partial}{\partial q}+(b+ic)(b'+ic')(\beta+i\gamma)^2e^{2i\varphi}e^{2qe^{i\varphi}}\frac{\partial}{\partial p}\otimes\frac{\partial}{\partial p}+$$
$$+[(db'-d'b)+i(dc'-d'c)](\beta+i\gamma)e^{i\varphi}e^{qe^{i\varphi}}\frac{\partial}{\partial p}\otimes\frac{\partial}{\partial q}$$
On the other hand,  $\langle F',[A,B]\rangle=[(db'-d'b)+i(dc'-d'c)](\beta+i\gamma)e^{i\varphi}e^{qe^{i\varphi}}$.\par  This implies (18). The other cases are proved similarly. 
\end{itemize}
 The proposition is hence completely proved.\hfill$\square$\par
\begin{definition}  Each chart $ ( \Omega_F ,   \psi^{-1})$  satisfying  1. and  2. of  the   Proposition $(3.1.1)$  is called an  adapted chart on $\Omega_F .$
\end{definition}
\subsection{Computation of operators $\hat \ell_A$ }
Since  $\widetilde A\circ \psi(p,q)=  \Phi(q).p+\Psi(q),$ for $A\in\mathfrak g$,  one can prove that :$$P^r(\tilde A,\tilde B)=0\quad  \forall r \geq 3,\forall A,B \in\mathfrak g.$$ From this   we have the following proposition
\begin{proposition}
With $A,B \in\mathfrak g$,  the Moyal $\star$-product satisfies the relation:
\begin{equation}i\tilde A\star i\tilde B - i\tilde B\star i\tilde A=i\widetilde{[A,B]}
\end{equation}
\end{proposition}
For each $A\in \mathfrak g$ and the corresponding Hamiltonian function $\tilde A$, we denote by  $ \ell_A$   the operator  acting on dense  space $L^2(\mathbb R^2,\frac{dpdq}{2\pi})$ of smooth function by left $\star$-multiplication  by $i\tilde A\quad $, i.e  $\ell_A(f)=i\tilde A\star f $. The relation in the Proposition (3.2.1) gives us 
\begin{corollary} 
\begin{equation}
\ell_{[A,B]}=\ell_A\star\ell_B - \ell_B\star\ell_A:=[\ell_A,\ell_B]^{\star}
\end{equation}
\end{corollary}
This implies that the corresponding $A\in\mathfrak g \mapsto\ell_A=i\tilde A\star.$ is a representation of the Lie algebra $\mathfrak g$ on the space $C^{\infty}(\Omega_F)[[\frac i 2]]$ of formal power series.\par
Let us denote by $\mathcal F_p(f)$ the partial Fourier transform (  define on $\L^2(\mathbb R^2,dpdq/2\pi)$, for example  ) of the function $f$ from the variable $p$ to the variable $x$ (see e.g   \cite{meisevogt}), i.e
$$\mathcal F_p(f)(x,q) := \frac{1}{\sqrt{2\pi}}\int_{\mathbb R} e^{-ipx} f(p,q)dp $$ and  $ \mathcal F^{-1}_p(f) (p,q)$ the inverse Fourier transform.\par 
We have  following obvious identities:
 \begin{lemma} With  $f,\varphi \in \L^2(\mathbb R^2,dpdq/2\pi)$ 
\begin{itemize}
\item[i.] $\partial_p \mathcal F^{-1}_p(f) = i \mathcal F^{-1}_p(x.f)$
\item[ii.]  $\mathcal F_p(p.\varphi) = i \partial_x\mathcal F_p(\varphi )  $
\item[iii.]  $P^r(\tilde A,\mathcal F^{-1}_p(f))=(-1)^r\partial_q^r(\Psi)\partial_p^r\mathcal F^{-1}_p(f) $  \quad $\forall r\geq 2.$
\end{itemize}
\end{lemma}
Now we denote $\hat\ell_A:=\mathcal F_p\circ\ell_A\circ\mathcal F^{-1}_p$ with $ A\in\mathfrak g.$
\begin{definition} Let    $  \Omega_F $  be K-orbit of co-adjoint representations of Lie group  $\G $. With  $A$ running over the Lie algebra $\mathfrak g=\Lie(\G),$ $(\Omega_F,\hat\ell_A)$  is called quantum co-adjoint orbits  of  Lie group $\G$.  
\end{definition}

\begin {theorem}
For each $A \in \mathfrak g$ and for each compactly supported $C^\infty$-function $f \in C^\infty_c(\mathbb R^2)$, we have 
$$\hat{\ell}_A(f) =\Phi(q-\frac x2)(\frac 12 \partial_q - \partial_x)f+i\Psi(q-\frac x2)f $$ and setting new variables $s=q-\frac x2,  t=q+\frac x2$, then
\begin{equation}
\hat \ell_A(f)=\Phi(s)\frac {\partial f}{\partial s}+i\Psi(s)f \vert_{(s,t)}\quad\mbox{i.e:}\quad\hat \ell_A=[\Phi(s)\frac {\partial }{\partial s}+i\Psi(s)] \vert_{(s,t)}
\end{equation}
\end{theorem}
{\it Proof.}   It is easy to see that:

$$P^0(\tilde A,\mathcal F^{-1}_p(f))=[\Phi(q).p+\Psi(q)]\mathcal F^{-1}_p(f)$$
$$P^1(\tilde A,\mathcal F^{-1}_p(f))=\Phi(q)\partial_q\mathcal F^{-1}_p(f) -[p\partial_q\Phi(q)+\partial_q\Psi(q)]\mathcal F^{-1}_p(f)$$
$$P^r(\tilde A,\mathcal F^{-1}_p(f))=(-1)^r\partial_q^r\Psi\partial_p^r\mathcal F^{-1}_p(f)   \quad \forall r\geq 2.$$

From this and  lemma (3.2.3), we have:$$\hat{\ell}_A(f) = \mathcal F_p \circ \ell_A \circ \mathcal F^{-1}_p(f)=i\mathcal F_p(\tilde A\star\mathcal F^{-1}_p(f))=i\mathcal F_p\left(\sum_{r \geq 0} \left(\frac{1}{2i}\right)^r\frac {1}{r!} P^r(\tilde A,\mathcal F^{-1}_p(f))\right)=$$
$$=i\mathcal F_p\{[\Phi(q).p+\Psi(q)]\mathcal F^{-1}_p(f)+\frac{1}{1!}\frac{1}{2i}(\Phi(q)\partial_q\mathcal F^{-1}_p(f) -[p\partial_q\Phi(q)+\partial_q\Psi(q)]\mathcal F^{-1}_p(f))+\dots+$$ $$+\frac{1}{r!}(\frac{1}{2i})^r(-1)^r\partial_q^r\Psi\partial_p^r\mathcal F^{-1}_p(f))+\dots\}=\Phi(q-\frac x2)(\frac 12 \partial_q - \partial_x)f+i\Psi(q-\frac x2)f $$
Theorem is proveed.\hfill$\square$\par
 As a direct consequence  of  symbol $\hat\ell_A,$  we have 
\begin{corollary}
$\forall A,B \in\mathfrak g ,$
$$\hat\ell_A\circ\hat\ell_ B- \hat\ell_B\circ\hat\ell_ A =\hat\ell_{[A, B]}$$
\end{corollary}
From theorems 2.0.3 and 3.2.5  we obtained the quantum half planes ,the quantum planes,  the quantum hyperbolic cylinders,  quantum hyperbolic paraboloids \dots of the corresponding groups. At the same time,  we  have also  unitary representations of these groups :$$T(exp A)=\exp(\hat\ell_A)=\exp\Big([\Phi(s)\frac {\partial }{\partial s}+i\Psi(s)] \vert_{(s,t)}
\Big).$$

\section{The case of groups    $\G_{4,2,3(\frac{\pi}{2})},\G_{4,2,4},\G_{4,3,4(\frac{\pi}{2})},\G_{4,4,1}$.} 
\subsection{The local diffeomorphisms}
For the Lie group of affine transformations of the complex straight line \quad $\G_{4,2,4}=\Aff(\mathbb C)$\quad (see \cite{diephai2}), we replaced the global diffeomorphism $\psi$  by  a local diffeomorphism and  obtained
\begin{proposition}Fixing the local  diffeomorphism $\psi_k (k \in {\mathbb Z})$
$$\psi_k: {\mathbb C} \times {\mathbb H}_k \longrightarrow {\mathbb C} \times {\mathbb C}_k$$
                           $$(z,w) \longmapsto (z,e^w), $$
 we have     \begin{enumerate}
        \item For any element $A \in \aff(\mathbb C)$, the corresponding Hamiltonian function $\widetilde{A}$  in local coordinates $(z,w)$ of the orbit $\Omega_F$  is of the  form
$$\widetilde A\circ\psi_k(z,w) = \frac{1}{2} [\alpha z +\beta e^w + \overline{\alpha} \overline{z} + \overline{\beta}e^{\overline {w}}]$$
       \item In local coordinates $(z,w)$ of the orbit\quad $\Omega_F$, the  Kirillov form $\omega$ is of the  form $$\omega = \frac{1 }{2}[dz \wedge dw+d\overline {z} \wedge d\overline {w}].$$
    \end{enumerate}
\end{proposition}
Analogously,  for the groups  $\G_{4,2,3(\varphi)},\G_{4,3,4(\varphi)}$ with $\varphi=\frac{\pi}{2}$,  we also replace the global diffeomorphism $\psi$  by  a local diffeomorphism $\psi_k (k \in {\mathbb Z}).$ \par
\begin{itemize}
\item   For $\G=\G_{4,2,3(\frac{\pi}{2})};\quad  \Omega_F=\{(\alpha,(\beta+i\gamma)e^{is},t)\vert s,t \in \mathbb R \}$,
$$\psi_k : \mathbb R\times (2k\pi,2\pi+2k\pi) \longmapsto \Omega_F .$$ 
$$(p,q)\longmapsto (\alpha,(\beta+i\gamma)e^{iq},p)$$
Then  the corresponding Hamiltonian fuction is $\tilde A\circ \psi_k(p,q)=dp+(b+ic)(\beta+i\gamma)e^{iq}+a\alpha. $\par
\item   For $\G=G_{4,3,4(\frac{\pi}{2})} ;\quad \Omega_F=\{((\alpha+i\beta)e^{is},\gamma e^{\lambda s},t)\vert s,t \in\mathbb R\}$,
$$\psi_k  :  \mathbb R\times (2k\pi,2\pi+2k\pi) \longmapsto \Omega_F .$$ 
$$(p,q)\longmapsto ((\alpha+i\beta)e^{iq},\gamma e^{q\lambda},p)$$
We have  $\tilde A\circ \psi_k(p,q)=dp+(a+ib)(\alpha+i\beta)e^{iq}+c\gamma e^{q\lambda}.$\par
\end{itemize}
At last,  for the Lie group $\G_{4,4,1}=\mathbb R\ltimes_j\mathbb H_3,$ which is not  exponential group, we obtain
\begin{proposition}
Let us denote $\mathbf I_k=(2k\pi,2\pi+2k\pi),k\in \mathbb Z$.
 Each non-trivial orbit $\Omega_F$ (in $\mathfrak g^*$) of co-adjoint representation of\quad $G_{4,4,1}$  admits local charts  $(\mathbb R\times\mathbf I_k,\psi^{-1}_k)$  or  $(\mathbb R_{\pm}\times\mathbf I_k,\psi^{-1}_k)$  such that:
\begin{enumerate} 
     \item If  $\Omega_F$ is defined by  (11)  and  $A\in \mathfrak g_{4,4,1},$  then $$\tilde A\circ\psi_k(p,q) =dp+\frac 1 2[a(\alpha+i\beta)+b(\beta-i\alpha)]e^{iq}+\frac 1 2[a(\alpha-i\beta)+b(\beta+i\alpha)]e^{-iq}$$
\begin{equation} =dp+(a\alpha+b\beta)cos q+(b\alpha-a\beta)sin q \end{equation}
and the Kirillov  form  then is  $\omega=dp\wedge dq.$
      \item  If  $\Omega_F$ is  defined by  (12) and  $A\in \mathfrak g_{4,4,1},$  then 
$$\tilde A\circ\psi_k(p,q)=\frac{d}{2\gamma}p^2+[\frac {a}{2}(e^{iq}+e^{-iq})+\frac{b}{2i}(e^{iq}-e^{-iq})]p+c\gamma+d\delta -d\frac {\alpha^2+\beta^2}{2\gamma}$$
\begin{equation}=\frac{d}{2\gamma}p^2+(acos q+bsin q)p+c\gamma+d\delta -d\frac {\alpha^2+\beta^2}{2\gamma}\end{equation}  and  we have the Kirillov form to be   $\omega=\frac{\gamma}{p}dp\wedge dq.$  
\end{enumerate}
\end{proposition}
{\it Proof.}
{\bf 1.} We consider  the  following  diffeomorphism: 
$$ \psi_k:\qquad \mathbb R \times \mathbf I_k\longrightarrow \Omega_F $$
$$(p,q)\mapsto(\frac 1 2(\alpha+i\beta)e^{iq}+\frac 1 2(\alpha-i\beta)e^{-iq};\frac 1 2(\beta-i\alpha)e^{iq}+\frac 1 2(\beta+i\alpha)e^{-iq};0;p)$$
With each $F'\in \Omega_F$,$$F'=[\frac 1 2(\alpha+i\beta)e^{iq}+\frac 1 2(\alpha-i\beta)e^{-iq}]X^*+[\frac 1 2(\beta-i\alpha)e^{iq}+\frac 1 2(\beta+i\alpha)e^{-iq}]Y^*+pT^*$$ and $ A=aX+bY+cZ+dT \in \mathfrak g_{4,4,1},$ we have $$\tilde A(F')=\langle F',A\rangle=\frac a 2[(\alpha+i\beta)e^{iq}+(\alpha-i\beta)e^{-iq}]+\frac b 2[(\beta-i\alpha)e^{iq}+(\beta+i\alpha)e^{-iq}]+dp.$$
It  follows that $$\xi_A(f)=d\frac{\partial f}{\partial q}-\frac i 2[a(\alpha+i\beta)+b(\beta-i\alpha)]e^{iq}+\frac i 2[a(\alpha-i\beta)+b(\beta+i\alpha)]e^{-iq}\frac{\partial f}{\partial p}.$$
By analogy,$$\xi_B(f)=d'\frac{\partial f}{\partial q}-\frac i 2[a'(\alpha+i\beta)+b'(\beta-i\alpha)]e^{iq}+\frac i 2[a'(\alpha-i\beta)+b'(\beta+i\alpha)]e^{-iq}\frac{\partial f}{\partial p}.$$                          Thus,
\begin{equation}\xi_A\otimes\xi_B=dd'\frac{\partial}{\partial q}\otimes\frac{\partial}{\partial q}-\frac 1 4\{[a'(\alpha+i\beta)+b'(\beta-i\alpha)]e^{iq}-[a'(\alpha-i\beta)+b'(\beta+i\alpha)]e^{-iq}\}\end{equation}$$\times\{[a(\alpha+i\beta)+b(\beta-i\alpha)]e^{iq}-[a(\alpha-i\beta)+b(\beta+i\alpha)]
e^{-iq}\}\frac{\partial}{\partial p}\otimes\frac{\partial}{\partial p}+\frac 1 2\{[(db'-d'b)(\alpha+i\beta)+$$   $$+(ad'-a'd)(\beta-i\alpha)]e^{iq}+[(db'-d'b)(\alpha-i\beta)+(ad'-a'd)(\beta+i\alpha)] e^{-iq}\}\frac{\partial}{\partial p}\otimes\frac{\partial}{\partial q}$$
On the other hand,    $[A,B]=(db'-d'b)X+(ad'-a'd)Y+(ab'-a'b)Z $  implies 
\begin{equation}\langle F',[A,B] \rangle
=\frac 1 2[(db'-d'b)(\alpha+i\beta)+(ad'-a'd)(\beta-i\alpha)]e^{iq}+\end{equation}$$ +\frac 1 2[(db'-d'b)(\alpha-i\beta)+(ad'-a'd)(\beta+i\alpha)]e^{-iq}
=\frac 1 2\{[(db'-d'b)(\alpha+i\beta)+$$ $$+(ad'-a'd)(\beta-i\alpha)]e^{iq}+[(db'-d'b)(\alpha-i\beta)+(ad'-a'd)(\beta+i\alpha)] e^{-iq}$$
(37) and (38) imply that  the Kirillov  form is  $\omega=dp\wedge dq$.\par
{\bf 2.}For  the  case $\gamma \ne 0$ we chose 
$$\psi_k:\qquad\qquad (\mathbb R_+\times \mathbf I_k  )\longrightarrow \Omega_F $$
  $$(p,q)\longmapsto (pcos q,psin q,\gamma,\frac {1}{2\gamma}(p^2+2\gamma\delta-\alpha^2-\beta^2))$$
$$\left(\mbox{or}\quad  \psi_k:\qquad\qquad (\mathbb R_-\times \mathbf I_k)\longrightarrow \Omega_F   \right).$$
Then,    $\forall F'\in \Omega_F  ,\quad F'=pcos q X^*+psin q Y^*+\gamma Z^*+\frac {1}{2\gamma}(p^2+2\gamma\delta-\alpha^2-\beta^2)T^*$    and $A=aX+bY+cZ+dT\in   \mathfrak g_{4,4,1}$,   we have $$\tilde A\circ\psi_{(U,k)}=apcos q +bpsin q +c\gamma + \frac {d}{2\gamma}(p^2+2\gamma\delta -\alpha^2 -\beta^2)=$$
$$=\frac{d}{2\gamma}p^2+[\frac {a}{2}(e^{iq}+e^{-iq})+\frac{b}{2i}(e^{iq}-e^{-iq})]p+c\gamma+d\delta -d\frac {\alpha^2+\beta^2}{2\gamma}.$$
It  follows that 
$$\xi_A(f)=[\frac{d}{\gamma}p+\frac{1}{2i}((ai+b)e^{iq}+(ai-b)e^{-iq})]\frac{\partial f}{\partial q}-\frac{p}{2}[(ai+b)e^{iq}-(ai-b)e^{-iq}]\frac{\partial f}{\partial p}.$$
By analogy,
$$\xi_B(f)=[\frac{d'}{\gamma}p+\frac{1}{2i}((a'i+b')e^{iq}+(a'i-b')e^{-iq})]\frac{\partial f}{\partial q}-\frac{p}{2}[(a'i+b')e^{iq}-(a'i-b')e^{-iq}]\frac{\partial f}{\partial p}.$$ 
Thus, 
\begin{equation}
\xi_A\otimes\xi_B=\frac{p}{2}[(ai+b)e^{iq}-(ai-b)e^{-iq}\frac{p}{2}[(a'i+b')e^{iq}-(a'i-b')e^{-iq}]\frac{\partial }{\partial p}\otimes\frac{\partial}{\partial p}+\end{equation}$$+[\frac{d}{\gamma}p+\frac{1}{2i}((ai+b)e^{iq}+(ai-b)e^{-iq})][\frac{d'}{\gamma}p+\frac{1}{2i}((a'i+b')e^{iq}+(a'i-b')e^{-iq})]\frac{\partial }{\partial q}\otimes\frac{\partial}{\partial q}+$$
$$+\Biggl\{\frac{p^2}{2\gamma}\biggl[\Bigl((da'-d'a)i+(db'-d'b)\Bigr)e^{iq}-\Bigl((da'-d'a)i-(db'-d'b)\Bigr)e^{-iq}\biggr]+p(ab'-a'b)\Biggr\}\frac{\partial}{\partial p}\otimes\frac{\partial}{\partial q}.$$
On the other hands,
\begin{equation}\langle F',[A,B] \rangle
=(db'-d'b)p cos q+(ad'-a'd) psin q+(db'-d'b)\gamma =\end{equation}
$$=\frac{p}{2}\Bigl[\bigl((a'd-d'a)i+(db'-d'b)\bigr)e^{iq}-\bigl((da'-d'a)i-(db'-d'b)\bigr)e^{-iq}
\Bigl]+\gamma(ab'-a'b). $$
From (39) and (40) we see  that  (17)  is of  the form $\omega=\frac{\gamma}{p}dp\wedge dq.$\hfill$\square$\par
\subsection{Computation of operators $\hat \ell^{(k)}_A$}
It is easy  to prove that,
\begin{itemize}
\item If $\G=\G_{4,2,3(\frac{\pi}{2})} $ then
$$\hat\ell_A^{(k)}=\Big(d\frac{\partial}{\partial_s}+i(b+ic)(\beta+i\gamma)e^{is}+a\alpha\Big)\vert_{(s,t)}$$
\item If   $\G=\G_{4,3,4(\frac{\pi}{2})} $ then
$$\hat\ell_A^{(k)}=\Big(d\frac{\partial}{\partial_s}+i(a+ib)(\alpha+i\beta)e^{is}+c\gamma e^{s\lambda}\Big)\vert_{(s,t)}$$
\end{itemize}
In \cite{diephai2}, we  proved the following result for $\G_{4,2,4}$ \par
Let ${\mathcal F}_z$(f)   denote     the partial Fourier transform of the function $f$ from the variable $z=p_1+ip_2$ to the variable $\xi=\xi_1+i\xi_2$, i.e:
$${\mathcal F}_z(f )(\xi,w) = \frac{1 }{2\pi} \iint_{R^2} e^{-iRe(\xi \overline{z})} f(z,w)dp_1dp_2$$  and   $$\mathcal F_z^{-1}(f )(z,w) = \frac{1}{2\pi} \iint_{R^2} e^{iRe(\xi \overline{z})} f( \xi,w)d \xi_1d\xi_2$$   the inverse Fourier transform.
\begin{theorem}( see  \cite{diephai2},Proposition 3.4)\par
For each $A = \left(\begin{matrix}\alpha & \beta \cr 0 & 0 \cr\end{matrix}\right) \in \aff({\mathbb C}) $ and for each compactly supported $C^{\infty}$-function $f \in C_c^{\infty}({\mathbb C} \times {\mathbb H}_k)$, we have:
\begin{equation} \hat\ell_A^{(k)}(f) := \mathcal F_z \circ \ell_A^{(k)} \circ \mathcal F_z^{-1}(f) =\end{equation}$$= [\alpha (\frac{1}{2} \partial_w - \partial_{\overline \xi})f + \overline \alpha (\frac{1 }{2}\partial_{\overline w} - \partial_\xi)f + \frac{i}{2}(\beta e^{w-\frac{1}{2}\overline \xi} + \overline \beta e^{\overline w - \frac{1}{2} \xi})f] $$ i.e $$\hat {\ell}_A^{(k)} = \alpha\frac{ \partial }{\partial u}+ \overline{\alpha}\frac{ \partial  }{\partial{\overline{u}}}+ \frac{i }{2}(\beta e^{u}+\overline{\beta}e^{\overline{u}});  \quad  u = w - \frac{1}{ 2}\overline{\xi};  v = w + \frac{1 }{2}{\overline{\xi}}.$$
\end{theorem} 
Now we consider  group $\G_{4,4,1}$ :
\begin{theorem} For each $A\in\mathfrak g_{4,4,1}$ and for each compactly supported $C^{\infty}$-function $f\in C_0^{\infty}(\mathbb R\times\mathbf I_k)$,  we have:  \par 
\begin{enumerate}
\item If $\tilde A$ is defined by  (35 ) then 
$$\hat\ell_A^{(k)}(f)=
\Bigl(d\partial _sf+i[(a\alpha+b\beta)cos s+(b\alpha-a\beta)sin s]f\Bigr)\vert_{(s,t)}$$
\item If $\tilde A$ is defined by  (36 ) then 
$$\hat\ell_A^{(k)}(f)=i.\Bigl([\frac{d}{2\gamma}(i\partial_x+\frac{1}{2\gamma}\partial_x\partial_q)^2]f+\Gamma .f \Bigr)\vert_{(x,q)}-$$ $$- \frac 12\Bigl(\partial_x[(a-bi)\Delta(f)+(a+bi)\Delta^{-1}(f)] + \frac{i}{2\gamma}\partial_x[(a-bi)\Theta(f)+(a+bi)\Theta^{-1}(f)] \Bigr)\vert_{(x,q)} $$
where
$$\Gamma=c\gamma+d\delta-d\frac{\alpha^2+\beta^2}{2\gamma}$$
$$\Delta(f)=\exp\bigl[iq+\partial_x((\frac{x}{2\gamma}).f)\bigr]=e^{iq}.\sum_{r=0}\frac{1}{r!}\frac{\partial^r}{\partial_{x^r}}\Bigg((\frac{x}{2\gamma})^r.f\Bigg)$$          $$\Theta(f)=\exp\bigl[iq+\partial_x((\frac{x}{2\gamma}).\partial_qf)\bigr]=e^{iq}.\sum_{r=0}\frac{1}{r!}\frac{\partial^r}{\partial_{x^r}}\Bigg((\frac{x}{2\gamma})^r.\partial_qf\Bigg).$$
\end{enumerate}
\end{theorem}
To prove the theorem, we need the following obvious lemma,   which is a direct  conseqnence of the definition of \quad$ \mathcal F_p$  and $\mathcal F_p^{-1}$.
\begin{lemma} With     $\quad\forall  r\geq 1$
   \begin{itemize}
    \item$ \partial^r_{p^r}\Bigl(\mathcal F_p^{-1}(f)\Bigr)=i^r\mathcal F_p^{-1}(x^r.f)$
    \item$ \mathcal F_p\Bigl(p^r\mathcal F_p^{-1}(f)\Bigr)=i^r\partial^r_{x^r}(f) $
    \item $ \mathcal F_p\Bigl(p^r\partial^{r-1}_{p^{r-1}}\mathcal F_p^{-1}(f)\Bigr)=i^{2r-1}\partial^r_{x^r}(x^{r-1}.f).$
\end{itemize}\hfill$\square$\par

\end{lemma}
 
{\it Proof.} \par
First case is  proved like as Theorem 3.2.5.
  We  prove only the second case.\par
 One  can write (36) as 
$$\tilde A\circ\psi(p,q)=\frac{d}{2\gamma}p^2+\frac{p}{2}[(a-bi)e^{iq}+(a+bi)e^{-iq}]+c\gamma+d\delta-d\frac{\alpha^2+\beta^2}{2\gamma}$$ and remark that
  in the  coordinates $(p,q)$ correspond to the form $\omega=\frac{\gamma}{p}dp\wedge dq$,$$\Lambda^{-1} = \left(\begin{array}{cc} 0 & \frac{p}{\gamma}\\ -\frac{p}{\gamma} & 0 \end{array}\right).$$  Denoting  $\mathcal F^{-1}_p(f)=v,$ we have
$$P^0=\Bigl\{\frac{d}{2\gamma}p^2+\frac{p}{2}[(a-bi)e^{iq}+(a+bi)e^{-iq}]+c\gamma+d\delta-d\frac{\alpha^2+\beta^2}{2\gamma}\Bigr\}v$$
$$P^1=[\frac{dp^2}{\gamma^2}+\frac{p}{2\gamma}((a-bi)e^{iq}+(a+bi)e^{-iq})]\partial_qv-\frac{ip^2}{2\gamma}[(a-bi)e^{iq}-(a+bi)e^{-iq}]\partial_pv$$
$$P^2=\frac{1}{2\gamma^2}\{i^2[(a-bi)e^{iq}+(a+bi)e^{-iq}]p^3\partial^2_{p^2}v-2 i [(a-bi)e^{iq}-(a+bi)e^{-iq}]p^2\partial^2_{pq}v\}+
\frac{dp^2}{\gamma^3}\partial_{q^2}^2v$$
$$P^3=\frac{1}{2\gamma^3}\{{(-i)}^3[(a-bi)e^{iq}-(a+bi)e^{-iq}]p^4\partial^3_{p^3}v+3{(-i)}^2[(a-bi)e^{iq}+(a+bi)e^{-iq}]p^3\partial^3_{p^2q}v\}$$
$$P^4=\frac{1}{2\gamma^4}\{{(-i)}^4[(a-bi)e^{iq}+(a+bi)e^{-iq}]p^5\partial^4_{p^4}v+4 {(-i)}^3[(a-bi)e^{iq}-(a+bi)e^{-iq}]p^4\partial^4_{p^3q}v\}$$
By analogy, for  $\quad \forall r \geq 4$  we have 
$$P^r=\frac{1}{2\gamma^r}\Bigl
((-i)^r [(a-bi)e^{iq}+(-1)^r(a+bi)e^{-iq}]p^{r+1}\partial^r_{p^r}v\Bigr) +$$  
$$+\frac{r{(-i)}^{r-1}}{2\gamma^r}\Bigl
([(a-bi)e^{iq}+(-1)^{r-1}(a+bi)e^{-iq}]p^r\partial^r_{p^{r-1} q}v\Bigr).$$
As $\tilde {A}$  is defined by  (36), we obtain:
\begin{equation}
\hat\ell^{(k)}_A(f) := \mathcal F_p \circ \ell_A \circ \mathcal F^{-1}_p(f)  = \mathcal F_p(i\tilde{A}\star \mathcal F^{-1}_p(f))
= i\sum_{r \geq 0} \Bigl(\frac{1}{2i}\Bigr)^r\mathcal F_p\left( P^r(\tilde{A},\mathcal F^{-1}_p(f)\right)=\end{equation}
$$=i(c\gamma+d\delta-d\frac{\alpha^2+\beta^2}{2\gamma})\mathcal F_pv+i\Bigl[\mathcal F_p(\frac {d}{2\gamma}p^2 v)+\frac{1}{1!}\frac{1}{2i}
\mathcal F_p(\frac{d}{\gamma^2}p^2\partial_qv)+
\frac{1}{2!}{\Bigl(\frac{1}{2i}\Bigr)}^2\mathcal F_p(\frac{d}{\gamma^3}p^2\partial^2_{q^2}v)\Bigr]+$$
$$+\sum_{r \geq 0}\frac{1}{r!}{\Bigl(\frac{1}{2i}\Bigr)}^r\frac{1}{2\gamma^r}\Bigl\{(-i)^r[(a-bi)e^{iq}+(-1)^r(a+bi)e^{-iq}]\mathcal F_p(p^{r+1}\partial^r_{p^r}v)+r(-i)^{r-1}
[(a-bi)e^{iq}+$$ $$+(-1)^{r-1}(a+bi)e^{-iq}]\mathcal F_p(p^r\partial^r_{p^{r-1}q}v)\Bigr\}= 
i(c\gamma+d\delta-d\frac{\alpha^2+\beta^2}{2\gamma})f+i\frac{d}{2\gamma}(i\partial_x+\frac{1}{2\gamma}\partial_x\partial_q)^2f
+$$ $$+\frac{i}{2}
\sum_{r \geq 0}\frac{1}{r!}\Bigl(\frac{-1}{2\gamma}\Bigr)^r(a-bi)e^{iq}\mathcal F_p(p^{r+1}\partial^r_{p^r}v)
+\frac{i}{2}\sum_{r \geq 0}\frac{1}{r!}\Bigl(\frac{1}{2\gamma}\Bigr)^r(a+bi)e^{-iq}\mathcal F_p(p^{r+1}\partial^r_{p^r}v)+$$ $$+\frac{i}{2}\frac{1}{\gamma}
\sum_{r \geq 1}\frac{1}{(r-1)!}\Bigl(\frac{-1}{2\gamma}\Bigr)^{r-1}\frac{1}{2i}(a-bi)e^{iq}\mathcal F_p(p^r\partial^r_{p^{(r-1)q}}v)+
\frac{i}{2}\frac{1}{\gamma}\sum_{r \geq 1}\frac{1}{(r-1)!}\Bigl(\frac{1}{2\gamma}\Bigr)^{r-1}\frac{1}{2i}(a+bi)$$ $$\times e^{-iq}
\mathcal F_p(p^r\partial^r_{p^{(r-1)q}} v ) =i[c\gamma+d\delta-d\frac{\alpha^2+\beta^2}{2\gamma}
+i\frac{d}{2\gamma}(i\partial_x+\frac{1}{2\gamma}\partial_x\partial_q)^2]f+\frac{i}{2}(a-bi) e^{iq}$$
$$\times\sum_{r \geq 0}\frac{1}{r!}\Bigl(\frac{-1}{2\gamma}\Bigr)^r i^{2r+1}\partial^{r+1}_{x^{r+1}}(x^r.f)+
\frac{i}{2}(a+bi)e^{-iq}\sum_{r \geq 0}\frac{1}{r!}\Bigl(\frac{1}{2\gamma}\Bigr)^r i^{2r+1}\partial^{r+1}_{x^{r+1}}(x^r.f)+$$ $$+\frac{1}{4\gamma}(a-bi)e^{iq}
\sum_{r \geq 1}\frac{1}{(r-1)!}\Bigl(\frac{-1}{2\gamma}\Bigr)^{r-1}i^{2r-1}
\partial^{r-1}_{p^{r-1}}(x^{r-1}.\partial_qf)+$$ $$+\frac{1}{4\gamma}(a+bi)e^{-iq}\sum_{r \geq 1}\frac{1}{(r-1)!}\Bigl(\frac{1}{2\gamma}\Bigr)^{r-1}i^{2r-1}\partial^{r-1}_{p^{r-1}}(x^{r-1}.\partial_qf)=$$
$$=i[c\gamma+d\delta-d\frac{\alpha^2+\beta^2}{2\gamma}+i\frac{d}{2\gamma}(i\partial_x+\frac{1}{2\gamma}\partial_x\partial_q)^2]f
-\frac 12\partial_x\Bigl[(a-bi)e^{iq}.\sum_{r \geq 0}\frac{1}{r!}\Bigl(\frac{1}{2\gamma}\Bigr)^r\partial^r_{x^r}(x^r.f)+$$
$$+(a+bi)e^{-iq}.\sum_{r \geq 0}\frac{1}{r!}\Bigl(\frac{-1}{2\gamma}\Bigr)^r \partial^r_{x^r}(x^r.f)
+\frac{i}{2\gamma}(a-bi)e^{iq}.\sum_{r \geq 1}\frac{1}{(r-1)!}\Bigl(\frac{1}{2\gamma}\Bigr)^{r-1}\partial^{r-1}_{x^{r-1}}(x^{r-1}.\partial_qf)+$$
$$+\frac{i}{2\gamma}(a+bi)e^{-iq}.\sum_{r \geq 1}\frac{1}{(r-1)!}\Bigl(\frac{-1}{2\gamma}\Bigr)^{r-1}\partial^{r-1}_{x^{r-1}}(x^{r-1}.\partial_qf)\Bigr]=$$
$$=i[c\gamma+d\delta-d\frac{\alpha^2+\beta^2}{2\gamma}+i\frac{d}{2\gamma}(i\partial_x+\frac{1}{2\gamma}\partial_x\partial_q)^2]f- $$  $$-\frac 12\partial_x
\Bigl((a-bi)e^{iq}e^{\partial_x((\frac{x}{2\gamma}).f)} +
(a+bi)e^{-iq}e^{\partial_x((\frac{-x}{2\gamma}).f)}\Bigr)-$$ $$-
\frac{i}{4\gamma}\partial_x\Bigl((a-bi)e^{iq}e^{\partial_x
((\frac{x}{2\gamma}).\partial_qf)}+
(a+bi)e^{iq}e^{\partial_x((\frac{-x}{2\gamma}).\partial_qf)}\Bigr)=$$ $$
=i[c\gamma+d\delta-d\frac{\alpha^2+\beta^2}{2\gamma}+\frac{d}{2\gamma}(i\partial_x+\frac{1}{2\gamma}\partial_x\partial_q)^2]f-$$
$$-\frac 12\partial_x\Biggl ((a-bi)\exp\bigl[iq+\partial_x((\frac{x}{2\gamma}).f)\bigr]+(a+bi)\exp\bigl[-iq+\partial_x((\frac{-x}{2\gamma}).f)\bigr]\Biggr)-$$
$$-\frac{i}{4\gamma}\partial_x\Biggl((a-bi)\exp\bigl[iq+\partial_x((\frac{x}{2\gamma}).\partial_qf)\bigr]+(a+bi)\exp\bigl[-iq+\partial_x((\frac{-x}{2\gamma}).\partial_qf)\bigr]\Biggr)$$
The theorem is therefore completely proved.\hfill$\square$\par

Thus, we  obtained all  the operators $\quad\hat\ell_A;\quad  \hat\ell^k_A$, which provides (global or local) representations of the $\MD_4$-algebras.
At last,   as   $\hat\ell_A, \quad \hat\ell^k_A$ are representations of the  $\MD_4$-algebras,  we have operators :    $\exp(\hat\ell_A);\quad \exp(\hat\ell^k_A)$
 are  representations of the corresponding connected and simply connected  $\MD_4$-groups. We say that they are the representations of $\MD_4$-groups   arising from the reduction of the procedure of deformation quantization.\par
\vfil\eject

\centerline{ ACKNOWLEDGMENT}
The author would like to acknowledge the support of Seminar "Deformation Quantization and Applications" under the supervision by Professor Do Ngoc Diep and would like to thank Dr. Nguyen Viet Dung for his useful comments .

\end{document}